\magnification=\magstep0 

\hsize=15.7truecm

\voffset=.8truecm          

\hoffset=.3truecm \baselineskip=1.2\normalbaselineskip

\overfullrule=0pt       

\def\r{\mathop{{\rm I}\mskip-4.1mu{\rm R}}\nolimits}
    
\def\C{\mathop{{\rm I}\mskip-9.0mu{\rm C}}\nolimits}

\def\non{{\vskip.5cm\noindent}}

\def\qed{\hfill $\sqcap \hskip-6.5pt \sqcup$\bigbreak}

\centerline{\bf ON THE COMPLEX ZEROS OF THE RIEMANN ZETA-FUNCTION}

\quad

\centerline{G. PUGLISI}

\quad

\quad

\centerline{\bf Abstract}

{\it\quad The purpose of this paper is to prove that the so-called Quasi-Riemann Hypothesis for}

{\it\quad the Zeta-function\ implies the Riemann Hypothesis .}

\quad

\quad

\centerline{\bf Introduction}

\quad

We shall be concerned with the proof of the following

\quad

{\bf Theorem.}\ \ {\it Let $\ \rho=\beta+i\gamma\ $ be the complex zeros of the Riemann zeta function $\ \!\zeta(s)\ . $}

\qquad\qquad\quad\ {\it If}
$$\sup_{\rho}\beta<1$$

\qquad\qquad\quad\ {\it then}
$$\beta = {1\over2}$$

\qquad\qquad\quad\ {\it for any $\ \!\rho\ .$} 

\quad

Before proceeding to prove this statement, I want to give here a detailed exposition of my

method.

Put
$$M_V(s)\ =\ \sum_{n\le V}{\mu(n)\over n^s}\quad\ \big(s=\sigma+it\big)\ .$$

The first object is to define a class of analytic functions $F_V(s)$ such that, when $V$ is large 
\item{(a)} $F_V(s)$ are real for real $s\ \!.$

\item{(b)} $F_V(\rho)=0\ $ for any complex zero $\ \rho=\beta+i\gamma\ .$

\item{(c)} $F_V(s)\ $ are regular functions and converge to 1 as $\ \!V\rightarrow+\infty\ $ uniformly with respect to$\ s\ $

in any bounded region on the right of the line $\ \sigma=\sup_{\rho}\beta$

Using Rouch\'e's theorem, it is not difficult to prove that $\ \!M_V(s_V) = 0\ $for a suitable $\ \!s_V\in$

$\r\ \!$ with $\ |s_V-1|\le V^{a-1}\ \!,\ \sup_{\rho}\beta < a < 1\ \!$ and$\ V\ $sufficiently large. This is performed by

lemma 2 and leads to definition (10), i.e.
$$F_V(s)\ \!=\ \!\zeta(s) M_V(s+s_V-1)\ .$$

Thus (a) and (b) are verified, while lemma 3 states
$$F_V(s)\ \!=\ 1\ \!+\ \!O_{\epsilon}\Big(V^{a-\sigma-r+\epsilon}\ \!(1+|t|)^{\epsilon}\Big)\qquad\ \ \big(a<\sigma+r<1\ \ ,\ \ 0<r \le a-\sup_{\rho}\beta\big)\leqno(i)$$

so all these functions satisfy (c).

We next consider a zero $\ \!\rho_0=\beta_0+i\gamma_0\ \!$ with $\ \beta_0>\sup_{\rho}\beta - \epsilon\ ,\ \gamma_0>10^{10}\ $ and put
$$\ a-\beta_0+r = z_0\ \ \!.\leqno(ii)$$

Define further
$$G_{UV}(s)\ \!=\ \!U^{s-\beta_0-z_0}\big(F_V(s+i\gamma_0)+F_V(s-i\gamma_0)\big)/2\ \ \!.\leqno(iii)$$

Then, by $(i),\ (iii)$
$$\eqalign{&\qquad\qquad\qquad\qquad\qquad\qquad G_{UV}(\beta_0)\ \!=\ \!0\cr
&G_{UV}(\beta_0 +z_0)\ \!=\ \!1\ \!+\ \!O_{\epsilon}\Big(V^{a-\beta_0-z_0-r+\epsilon}\ \!{\gamma_0}^{\epsilon}\Big)\ \!=\ \!1\ \!+\ \!O_{\epsilon}\Big(V^{-2r+\epsilon}{\gamma_0}^{\epsilon}\Big)\ .}\leqno(iv)$$

Moreover, given any integer $\ \!J>1\ \!$ one may find a  number $\ \!z_*\ ,\ 0<z_*< z_0\ \!$ such that
$$\eqalign{&G_{UV}(\beta_0)\ \!=\ \!G_{UV}(\beta_0+z_0)\ \!+\ \!\sum_{j=1}^{J-1}{(-z_0)^j\over j!}\ \!D^j\ \!G_{UV}(\beta_0+z_0)\ \!+\ \!{(-z_0)^J\over J!}\ \!D^J\ \!G_{UV}(\beta_0+z_*)\ .}\leqno(v)$$

On the other hand, using classical tools it is not difficult to show that
$$\eqalign{&F_V(s)\ \ll_{\epsilon}\ \big(V(1+|t|)\big)^{a-\sigma-r+\epsilon}\qquad\qquad\qquad\quad\big(\sigma<\sup_{\rho}\beta\big)\ \ .}\leqno(vi)$$

This bound is proved in lemma 4. Therefore, by appealing to Cauchy's inequality, one may

deduce from $(iii),\ (vi)$
$$\eqalign{&\qquad\qquad\quad\ {D^j G_{UV}(\beta_0+z_0)\over j!}\ \!-\ \!{\log^jU\over j!}\ \!\ll_{\epsilon}\ \!{(V\gamma_0)^{\epsilon}U^r\over r^jV^r}\ \ ,\cr
&\ {D^J G_{UV}(\beta_0+z_*)\over J!}\ \!-\ \!{U^{z_*-z_0}\ \!\log^JU\over J!}\ \!\ll_{\epsilon}\ \!{\big(V\gamma_0/U\big)^{z_0-z_*} (V\gamma_0)^{\epsilon}\over\ r^J}\ \!\Big({U\over V}\Big)^r}\leqno(vii)$$

for $\ \!j\ge1\ \!$ and $\ \! r\ \!$ as in $\ \!(i),\ \!(ii)\ \!.\ \!$ But these estimates are of little use when $j$ increases
and $z_*$

is not close to $z_0\ \!.\ \!$ We then proceed as follows.

Put
$$F_V(s)\ =\ \!\sum_{n=1}^{+\infty}{c_n\over n^s}\quad\big(\sigma>1\big)\quad\ ,\quad\ p_j(u)\ \!=\ \!{1\over j!}\ e^{-u}\ \!u^j\quad \big(u\in\r\big)$$

so that, by $\ \!(iii)$
$${(-z_0)^j\over j!}D^j\ \!G_{UV}(s)\ \!=\ \!{U\over2}^{s-\beta_0-z_0}\sum_{n=1}^{+\infty}\Big({c_n\over n^{s-z_0+i\gamma_0}}+{c_n\over n^{s-z_0-i\gamma_0}}\Big)\ \!p_j\big(z_0\log(n/U)\big)\quad\ \ \big(\sigma>1\big)\ \ \!.$$

The main problem is the treatment of the sum
$${{z_0}^j\over j!}\sum_{U<n\le V^2}\Big({c_n\over n^{s-z_0+i\gamma_0}}+{c_n\over n^{s-z_0-i\gamma_0}}\Big)\ \!p_j\big(z_0\log(n/U)\big)\qquad \big(1\le j\le J\big)$$

if $\ \!s=\beta_0+z_0\ \!$ and if $\ \!s=\beta_0+z_*\ \!.\ \!$ The final result is \big(see lemma 5\big) 
$$\eqalign{&\sum_{1\le j<J}{(-z_0)^j\over j!}\ \!D^j G_{UV}(\beta_0+z_0)\ \!+\ \!{(-z_0)^J\over J!}\ \!D^J G_{UV}(\beta_0+z_*)\ =\ \!\sum_{j=1}^{J-1}{(-z_0\log U)^j\over j!}\ +\cr
&\qquad\qquad\qquad\ +\ \!{U^{z_*-z_0}(-z_0\log U)^J\over J!}\ \!+\ \!{\cal R}\big(J, U, V, \beta_0,\gamma_0, z_0, z_*\big)}\leqno(viii)$$

where, if$\ U\ $is suitably chosen in terms of$\ V,\ $then$\ {\cal R}\ $is small when$\ V\ $and$\ J\ $are as large as

we need.

Furthermore, when $\ \!J\ge 2z_0\log U + 2\ \!$ is even, lemma 1 shows that
$$\sum_{j=1}^{J-1}{\big(-(z_0\log U\big)^j\over j!}\ \le\ -{(z_0\log U)^{J-1}\over2(J-1)!}\ \ .\leqno(ix)$$

Suppose now $\ \!V\ge {\gamma_0}^{2/r}\ \!$ and take $\ \!U = V^{2/3}\ \!,\ J = 2[z_0\log U+2]\ \!.$

If $(vi) ,\ \!$  holds for a $z_*\in \big(0, z_0\big)\ \!,$ then$\ (iv),\ \!(viii)\ $give
$$\sum_{j=1}^{J-1}{(-z_0\log U)^j\over j!}\ \!+ \ \!{U^{z_*-z_0}(-z_0\log U)^J\over J!}\ \!=\ \!-1\ \!+\ \!o(1)$$

as $\ \!V\ \!\rightarrow\infty\ \!.$ Therefore, by$\ (ix)\ $and the second statement of lemma 1
$$3\ \!>\ \!{(z_0\log U)^{J-1}\over(J-1)!}\ \!\Big(1-{2z_0\log U\over J}\Big)\ \!>\ \! {1\over e(J-1)}\Big({ez_0\log U\over J-1}\Big)^{J-1}\Big(1-{2z_0\log U\over J}\Big)$$

when $V$ is large enough. But this is impossible, since $\ \!2\le J-2z_0\log U\le 4\ .$

Then $\ \sup_{\rho}\beta\ = 1/2\ .$

\non

\centerline{\bf Proof of the theorem}

\non

We put
$$\ 0\ \!\le \sup_{\gamma> 10^3}\Big(\beta-{1\over2}\Big)\ \!=\  b\ \!\le\ \! {1\over2}$$

and suppose
$$ 0<b<{1\over2}\quad\ .\leqno(1)$$

then 
$$\eqalign{&(*)\quad either\quad\ {1\over2} > b >b'=\limsup_{\gamma\rightarrow+\infty}\Big(\beta-{1\over2}\Big)\ge 0\cr
&(**)\ or\qquad\quad\ b = b'\ \ .}$$

Moreover
$$\eqalign{&if\ \ b\ \ is \ as \ in\ \ (1) ,\ \!(*)\quad\ take\quad\rho_0 = b+{1\over2}+i\gamma_0\quad\ where\cr
&\qquad\qquad\qquad\qquad\quad\ \gamma_0 = \max\Big\{\gamma>10^{10}\ \!|\ \exists\ \!\rho = b+{1\over2}+i\gamma\Big\}\cr
&if\ \ b\ \ is\ as\ in\ \ (1) ,\ \!(**)\quad then\quad\exists\ \!\rho_0 = \beta_0+i\gamma_0\ \ \!:\ \ \!\beta_0\ge b+(1-2\epsilon)/2\ \ ,\ \ \gamma_0>10^{10}\cr
&\qquad\qquad\qquad\qquad\qquad\qquad\ \!where\qquad 0<\epsilon\le 10^{-4}\ .}\leqno(2)$$

Following Bombieri \big(see [1], p. 46\big), define
$$\eqalign{&\qquad\qquad\qquad\qquad\quad\ p_j(u)\ =\ {1\over j!}\ e^{-u}\ \!u^j\quad \big(u\in\r\big)\cr
&so\ \ that\qquad\quad\ \sum_{j=0}^{\infty}\ \!p_j(u) = 1\quad\ ,\quad\ p'_j(u)\ =\ p_{j-1}(u) - p_j(u)\quad\ \ \big(j\ge 1\big)\ \ .}\leqno(3)$$

{\bf Lemma 1.}\ \ {\it If $\ \!J\ge2\ \!$ is an even integer and $\ u\le 0\ \!,\ \!$ then}
$$\sum_{j=1}^{J-1}p_j(u)\ \!\le {1\over2}\ \!p_{J-1}(u)\ .$$ 

\qquad\qquad\quad\ \ {\it Furthermore, if $\ \!j\ge 1$}
$$j!\ \!\le\ \!j^{j+1}\ \!\!e^{1-j}\ \ .$$

{\bf Proof.}\ \ The latter inequality is trivial for$\ j=1\ $while, if it is true for some$\ j\ \!,\ $then
$$(j+1)!\ \le\ \!(j+1)\ \!j^{j+1}e^{1-j}\ \!=\ e\Big({j\over j+1}\Big)^{j+1}\!(j+1)^{j+2}\ \!e^{-j}\ \!\le\ (j+1)^{j+2}\ \!e^{1-(j+1)}$$

since $\ \big(j/(j+1)\big)^{j+1}\ $ is an increasing sequence which has limit $1/e\ .$

As regards the former, which is true for$\ J=2\ $and$\ u\le0\ \!,\ $ it is sufficient to prove that

$$ {1\over2}\ \!p_{J-1}(u)+ p_J(u)\ \!\le\ \!-{1\over2}\ \!p_{J+1}(u)\ \ .$$

By (3) this last inequality is equivalent to
$$\ u^2+2(J+1)\ \!\!u+J(J+1)\ge 0$$

or to
$$|u+J+1|\ge\sqrt{J+1}$$

which actually holds, since
$$J/2\le J+1-\sqrt{J+1}\qquad\big(J\ge-1\big)\ \ .$$

\qed

According to (1), let now $\ a,\ r,\ v,\ s_0,\ T,\ \epsilon\ $ be real numbers such that

$$\eqalign{& 0<200\epsilon\le 2r\le\min(1-a\ \!,\ \!\ \ 1/50)\quad,\quad 1 > a\ge b+(1+2r)/2\quad,\quad 0<b<1/2\cr
&\qquad\qquad\qquad\quad w = a+iv\quad,\quad s_0 = a+r\quad,\quad2\le T\le v\le2T}\leqno(4)$$

and put
$$M_V(s)\ =\ \sum_{n\le V}{\mu(n)\over n^s}\qquad\ \big(V\ge2\big)\quad .\leqno(5)$$

\non

{\bf Lemma 2.}\ \ {\it  Let $\ \!\epsilon\ ,\ r\ ,\ a\ ,\ b\ \!$ be as in (4) and let $\ \!M_V(s)\ \!$ be defined by $(5)\ \!.\ \!$ If $\ V\ge V_0(\epsilon) ,$} 

\qquad\qquad\quad\ \ {\it then there exists a unique $\ s_V\in\big\{s\in\r : |s-1|< V^{a-1}\big\}\ \!$ such that}

$$\ M_V(s_V) = 0\ \ .$$

{\bf Proof.}\ \ Perron's formula (see [2] Lemma 3.12 and Lemma 3.19, with $\ \!a_n = \mu(n)\ \!,\ \psi(n)=$

$=\alpha = 1\ \!,\ x =V\ \!,\ T = W$\big) gives at once 
$$M_V(s)\ =\ {1\over2\pi i}\int_{c-iW}^{c+iW}\zeta^{-1}(s+z)\ {V^z\over z}\ dz\ +\ O\bigg({V^c\log V\over W}+{\log V\over V^{\sigma}}\bigg)\leqno(6)$$

where
$$c\ = \ \!\max(1-\sigma\ \!,\ 0) +(\log V)^{-1}\ \ ,\quad \sigma \ge a-r+\epsilon\quad,\quad W\ge2\ \ .$$

Also
$$\eqalign{&{1\over2\pi i}\int_{c-iW}^{c+iW}\zeta^{-1}(s+z)\ {V^z\over z}\ dz\ =\ \zeta^{-1}(s)\ +\ {1\over2\pi i}\Bigg(\int_{a-\sigma-r+{\epsilon\over2}-iW}^{a-\sigma-r+{\epsilon\over2}+iW}+\cr
&\qquad\ \ +\int_{a-\sigma-r+{\epsilon\over2}+iW}^{c+iW}-\int_{a-\sigma-r+{\epsilon\over2}-iW}^{c-iW}\Bigg)\ \zeta^{-1}(s+z)\ {V^z\over z}\ dz\ .}\leqno(7)$$

We now observe that the same argument which gives (14.2.6) in [2], assuming (1) instead of RH

leeds to
$$\zeta^{-1}(s)\ \!=\ \!O\big(|t|^{\epsilon}\big)\qquad\quad\big(\sigma>b+1/2\big)\ \ \!.$$

More precisely, when $\ \!\sup_{\rho}\beta = b+ 1/2\ \!$ the statement of Theorem 14.2 on pag 336 becomes
$$\log\zeta(s)\ \!=\ \!O\Big((\log t)^{{2-2\sigma\over1+2b}+\epsilon}\Big)$$

uniformly for $\ \!b + 1/2 <\sigma_0\le\sigma\le 1\ \!.\ \!$ This implies (14.2.5), (14.2.6) in the same range.

You must simply take the radii  $3/2 - b - \delta/2\ \!,\  3/2 - b- \delta\ \!$ in place of $\ \!3/2 - \delta/2\ \!,\  3/2 - \delta\ \!$ \big(see

(14.2.2)\big) and $C_3$ the circle with centre $\ \!\sigma_1 + it\ \!$  passing through the point  $\ \!b + 1/2  + \delta + it\ .$

Similarly, if, according to (1), (2), $1/2\ \!$ is replaced by $\ \!b + 1/2\ \!,\ \!$ then we have, by (4)

$$\ \ \!\zeta^{-1}(s+z)\ \ll_{\epsilon}\ (W+|t|)^{\epsilon/4}\qquad\big(|{\cal I}m
z|\le W\big)$$

uniformly for $\ \!{\cal R}e(z)+\sigma\ge a-r+\epsilon/2\ \!\ge\ \! b + (1+\epsilon)/2\ \!.\ \!$ Hence 
$$\eqalign{&\ \int_{a-\sigma-r+{\epsilon\over2}-iW}^{a-\sigma-r+{\epsilon\over2}+iW}\zeta^{-1}(s+z)\ {V^z\over z}\ dz\ \ll_{\epsilon}\ V^{a-\sigma-r+{\epsilon\over2}}\ \!(W+|t|)^{\epsilon/3}\cr
&\ \ \!\int_{a-\sigma-r+{\epsilon\over2}\pm iW}^{c\pm iW}\zeta^{-1}(s+z)\ {V^z\over z}\ dz\ \ll_{\epsilon}\ V^c W^{{\epsilon\over4}-1}(1+|t|)^{\epsilon/4}\ .}\leqno (8)$$

From (6), (7), (8) it follows that
$$M_V(s)\ \!=\ \!\zeta^{-1}(s)\ \!+\ \!O_{\epsilon}\Big(V^c\log V\ \!W^{-1}(W+|t|)^{\epsilon/4}+\ V^{a-\sigma-r+{\epsilon/2}}\ \!(W+|t|)^{\epsilon/3}\Big)\leqno(9)$$

where
$$W\ge 2\ \ ,\ \ \sigma\ge a-r+\epsilon\ \ ,\ \ c\ = \ \!\max(1-\sigma\ \!,\ 0) +(\log V)^{-1}\ .$$

Take $\ \ \!s\in\big\{s\in\C : |s-1|= V^{a-1}\big\}\ \ \!,\ \ \!W = V\ge V_0(\epsilon) \ge\big(6/\epsilon\big)^{1/(1-a)}\ \!.\ $ Then (4), (9) give, for 

a suitable $\ \!C(\epsilon)$
$$|M_V(s)-\zeta^{-1}(s)|\ \le\ C(\epsilon)\ \! V^{a-\sigma-r+5\epsilon/6}\ \!\le\ C(\epsilon)\ \! V^{a-1-r+\epsilon}$$

since $\ \!\sigma\ge 1-V^{a-1}\ge 1-\epsilon/6 \ge a+2r-\epsilon/6 > a-r+\epsilon\ ,\ \!$ while
$$|\zeta^{-1}(s)|\ \ge\ |s-1|/2\ =\ V^{a-1}/2\ >\  2C(\epsilon)V^{a-1-r+\epsilon}\ .$$

Therefore, by Rouch\'e's theorem, applied to the disk $\ \!\big\{s\in\C : |s-1|\le V^{a-1}\big\}$  
$$\sharp\big\{s\in\C : |s-1|<V^{a-1},\ \!M_V(s) = 0\big\}\ =\ \sharp\big\{s\in\C : |s-1|<V^{a-1},\ \!\zeta^{-1}(s) = 0\big\}\ =\ 1\ \ \!.$$

Hence $s_V$ is unique and then it is also real, since $\ \! M_V(s)\ \!$ is real for real $s$ (and consequently

has its non-real zeros in conjugate pairs).

This establishes lemma 2.
\qed

As a consequence of lemma 2, we define the following entire function 
$$\eqalign{&\quad\ \ F_V(s)\ =\  \zeta(s)\ \!M_V(s+s_V-1)\cr
&\ V\ ,\ M_V(s)\ ,\ s_V\ \ as\ \ in\ \ lemma\ \ 2\ \ .}\leqno(10)$$ 

We first prove

\non
{\bf Lemma 3.}\ \ {\it If$\ F_V(s)\ $is defined by (10) and $\ \!\epsilon\ ,\ r\ ,\ a\ \!$ are as in (4), then}
$$\eqalign{&\ \!F_V(s)\ \!=\ \!1\ \!+\ \!O_{\epsilon}\Big(V^{a-\min(1,\ \!\sigma+r)+\epsilon}\ \!(|t|+1)^{\epsilon}\Big)\cr
&\qquad\qquad for\qquad\ \sigma\ge a-r+2\epsilon\ \ .}$$
\non

{\bf Proof.} \ \ By lemma 2 $\quad\!|1-s_V|\ \le\ V^{a-1}\le\ \epsilon\quad$if $\ \ V\ge \epsilon^{1/(a-1)}\ .$

Then, when $\ \sigma\ge a-r+2\epsilon$
$$\min\big(\sigma,\ \sigma+s_V-1\big)\ \ge\ a+2\epsilon-r - |s_V -1|\ \ge\ a-r +\epsilon\ .$$ 

It then follows from (9) with $\ \!W=V\ \!$ that
$$\eqalign{&M_V(s) - M_V(s+s_V-1)\ \!=\ \!\zeta^{-1}(s) - \zeta^{-1}(s+s_V-1)\ \!+\ \!O_{\epsilon}\Big(\big(V^{a-\sigma-r+\epsilon}+V^{\epsilon-1}\big) (|t|+1)^{\epsilon/2}\Big)\ \!,\cr
&\qquad\ \ \big|\zeta^{-1}(s)-\zeta^{-1}(s+s_V-1)\big|\ \!=\ \!\bigg|\int_{s+s_V-1}^s{\zeta'\over\zeta^2}(z)\ dz\ \!\bigg|\ \!\le\ \!|s_V-1|\max_{{\cal R}e z\ge a-r+\epsilon\atop {\cal I}m z = t,\ \!|z-1|\ge \epsilon}\Big|{\zeta'\over\zeta^2}(z)\Big|\ .}$$

Also, by 14.2.5, 14.2.6\ \ in [2] \big(with $1/2$ replaced by $\ \!b + 1/2$\big) we have
$$\max\big(|\zeta(z)|\ \!,\ \!|\zeta(z)|^{-1} \big) \ll_{\epsilon} (|{\cal I}m z|+1)^{\epsilon/6}\qquad\big({\cal R}e z\ge a-r+\epsilon/2\ \ ,\ \ |z-1|\ge\epsilon/2\big)\leqno(11)$$

and applying Cauchy's inequality to $\zeta'(z)$ in the circle $\ \!|z-s|\le\epsilon/2\ \!$ \big(when $\ \!|s-1|\ge 2\epsilon$\big), we

obtain
$$\zeta'/\zeta^2(z)\ \ll_{\epsilon}\ \!(|t|+1)^{\epsilon/6}\ \!|\zeta^{-2}(z)|\ \!\ll_{\epsilon}\ \!({t}+1)^{\epsilon/2}\qquad\big({\cal R}e z\ge a-r+\epsilon\ \ ,\ \ |z-1|\ge\epsilon\big)\ .$$

Hence
$$\big|\zeta^{-1}(s)-\zeta^{-1}(s+s_V-1)\big|\ \ll_{\epsilon} V^{a-1}\big(|t|+1\big)^{\epsilon/2}$$

and then
$$\eqalign{&\big|F_V(s) - \zeta(s)\ \!M_V(s)\big|\ =\ |\zeta(s)|\ \!|M_V(s) - M_V(s+s_V-1)|\ \ll_{\epsilon}\cr
&\qquad\qquad\quad\ \ \ll_{\epsilon}\ (|t|+1)^{\epsilon}\Big(V^{a-\sigma-r+\epsilon} +\ \!V^{a-1}\Big)}\leqno(12)$$

when $\ \!|s-1|\ge 2\epsilon\ \!.\ $ Furthermore, again by (9) with $\ \!W=V\ $
$$\eqalign{(\zeta\ \!M_V)(s)\ \!=\ \!1 &+ O_{\epsilon}\Big(\big(V^{a-\sigma-r+\epsilon}+V^{\epsilon-1}\big) (|t|+1)^{\epsilon}\Big)\cr
&\quad\big(\sigma\ge a-r+\epsilon\big)\ \ .}\leqno(13)$$

The lemma now follows from (12), (13) if $\ \!|s-1|\ge2\epsilon\ \!.\ \!$ But it is also true when $\ \!|s-1|\le2\epsilon\ \!$ 

since $\ \!F_V(s) - 1\ \!$ is olomorphic.

\qed

\non

For $\ 2\le U\le V\ $ put, according to (10)
$$G_{UV}(s)\ \!=\ \!U^{s-s_0}\big(F_V(s+iv)+F_V(s-iv)\big)/2\qquad\qquad\big(s_0\ ,\ v\ \ as\ \ in\ \ (4)\big)\ \!.\leqno(14)$$ 

It is important to obtain sharp estimates for $\ \!D^jG_{U,V}(s) - U^{s-s_0}\log^jU\ \!$ when $\ \!s\ \!$ is near to $a\ \!.$

The following approach is based on a rather complicated argument involving the functions $p_j\ $

defined in (3). We begin by proving

\non

{\bf Lemma 4.}\ \ {\it Let $\ \epsilon\ ,\ r\ ,\ a\ ,\ w\ ,\ T\ $ be as in (4) and let $\ F_V(s)\ $ be defined by (10). If}
$$|\omega|\le a-1/2\quad,\quad 0\le{\cal R}e\omega\le a -1/2\quad,\quad{\cal R}e z\ge{\cal R}e\omega + 1/2 -a-r$$

\qquad\qquad\quad\ \ {\it then}
$$\max\Big(|F_V(w+r-\omega+ z)|\ ,\ \!|F_V(\overline w+r-\omega+ z)|\Big)\ \!\ll_{\epsilon}\ \!\big(V(T+|{\cal I}m z|)\big)^{\max(0,\ \!{\cal R}e (\omega-z)-2r)+4\epsilon}\ .$$

{\bf Proof.}\  According to (4), it suffices to prove the above inequality for $\ \!F_V(w+r-\omega+ z)\ \!.\ $

Suppose first that $\ \ \!-r +  2\epsilon\le {\cal R}e(z-\omega)+r\le 1-a+\log^{-1}\!V\ .\ $

Then, for $\ 2\le T_1\le V\ \!,\ \ V\ge V_0(\epsilon)\ \!,\ $ Perron's formula with $\ \!c=2-s_V-a +{\cal R}e(\omega -z)+r\ \!+$

$+2/(\log V)\ge\ \!1/\log V\ \!$ gives \big(see (5), (6), (7) (8)\big) 
$$\eqalign{&\ \!M_V(w+r+s_V-1-\omega+z)\ =\ {1\over2\pi i}\int_{2-s_V-a+{\cal R}e(\omega-z)-r+2\log^{-1}\!V-iT_1}^{2-s_V-a+{\cal R}e(\omega-z)-r+2\log^{-1}\!V+iT_1}\ \!{1\over\zeta}\big(w+r\ \!+\cr
&\qquad\quad\ \ \! +s_V-1-\omega+ z +\eta\big)\ V^{\eta}\ {d\eta\over\eta}\ \!+\ O\bigg({V^{1-a+{\cal R}e(\omega-z)-r}\log V\over T_1}\bigg)\ =\cr
&\quad\ \ =\ {1\over2\pi i}\Bigg(\int_{1-s_V+{\cal R}e(\omega-z)-2r+\epsilon-iT_1}^{1-s_V+{\cal R}e(\omega-z)-2r+\epsilon+iT_1} +\int_{1-s_V+{\cal R}e(\omega-z)-2r+\epsilon+iT_1}^{2-s_V-a+{\cal R}e(\omega-z)-r+2\log^{-1}\!V+iT_1}+\cr
&\ -\int_{1-s_V+{\cal R}e(\omega-z)-2r+\epsilon-iT_1}^{2-s_V-a+{\cal R}e(\omega-z)-r+2\log^{-1}\!V-iT_1}\Bigg)\ \!{1\over\zeta}(w+r+s_V-1-\omega+z+\eta)\ V^{\eta}\ {d\eta\over\eta}\ +\cr
&\ \!+\  {1\over\zeta}\big(w+r+s_V-1-\omega+ z\big) + O\bigg({V^{1-a+{\cal R}e(\omega-z)-r}\log V\over T_1}\bigg) \ll_{\epsilon}V^{{\cal R}e(\omega-z)-2r}\ \!\cdot\cr
&\quad\ \!\cdot\big((T+|{\cal I}m z|+T_1)\ \!V^2\big)^{\epsilon\over2} +\ \!\big(V^{1-a+{\cal R}e(\omega-z)-r}\log V\big)\ \!T_1^{-1} +\ \!\big(T+|{\cal I}m z|\big)^{\epsilon\over6} }$$

since, by lemma 2, $\ \!|1-s_V|\le 1/\log V\ \!$ and since $\ \!\zeta^{-1}(s)\ll_{\epsilon} (|t|+1)^{\epsilon/6}\ $ if $\ \sigma\ge a-r+\epsilon$

$\big($see (11)\big) . On the other hand, when $\ \!{1\over2}-a\le{\cal R}e(z-\omega)+r\le -r+2\epsilon\ \!$ we have
$$\eqalign{&\ \ \!M_V(w+r+s_V-1-\omega+z)\ =\ {1\over2\pi i}\Bigg(\int_{1-s_V+{\cal R}e(\omega-z)-2r+3\epsilon+iT_1}^{2-s_V-a+{\cal R}e(\omega-z)-r+2\log^{-1}\!V+iT_1}+\cr
&-\!\int_{1-s_V+{\cal R}e(\omega-z)-2r+3\epsilon-iT_1}^{2-s_V-a+{\cal R}e(\omega-z)-r+2\log^{-1}\!V-iT_1}\!+\int_{1-s_V+{\cal R}e(\omega-z)-2r+3\epsilon-iT_1}^{1-s_V+{\cal R}e(\omega-z)-r+iT_1}\Bigg)\ \!{1\over\zeta}(w+r\ \!+\cr
&\qquad+s_V -1-\omega+z+\eta)\ V^{\eta}\ {d\eta\over\eta}\ \!+\ \! O\bigg({V^{1-a+{\cal R}e(\omega-z)-r}\log V\over T_1}\bigg)\ \!\ll_{\epsilon}\cr
&\ \ll_{\epsilon}\ \! V^{{\cal R}e(\omega-z)-2r}\big((T+|{\cal I}m z|+T_1)V^6\big)^{\epsilon\over2}\ \!+\ \big(V^{1-a+{\cal R}e(\omega-z)-r}\log V\big)\ \!T_1^{-1}\ \!.}$$

Therefore, taking $\ T_1 = V^{1-a+r}\log V$ 
$$\eqalign{&M_V(w+r+s_V-1-\omega+z)\ \ll_{\epsilon}\ \!V^{\max(0,\ \!{\cal R}e(\omega-z)-2r)}\big(V^7(T+|{\cal I}m z|)\big)^{\epsilon\over2}\cr
&\qquad\qquad\quad\Big(\ \!{1\over2}-a\le{\cal R}e(z-\omega)+r\le 1-a+\log^{-1}\!V\ \!\Big)\ .}\leqno(15)$$

Also, by hypothesis, $\ \!a-r\ge (2b +1)/2\ge 1/2\ \!,\ \!$ so that $\ \!\ \!\zeta(it)\ \!\ll_{\epsilon}\ \!|t|^{{1\over2}+{\epsilon\over4}} \le\ \!|t|^{a-r+{\epsilon\over4}}\ \!.$

Then, since, by (11),
$$\ \zeta\big(a-r+\epsilon/2+it)\ \!\ll_{\epsilon}\ \!(1+|t|)^{\epsilon\over6}\quad\big(|t|\ge \epsilon\big)$$  

a well known convexity argument for the function $\ \!\mu(\sigma)\ \!$ (see [3] \S 5.1) gives
$$\eqalign{&\zeta(w+r-\omega+ z)\ \!\ll_{\epsilon}\big(T+|{\cal I}m z|\big)^{\max(0,\ \!{{\cal R}e(\omega-z)-2r\over 2(a-r)})+{\epsilon\over2}}\ \le\cr
&\qquad\qquad\qquad\quad\ \ \!\le\ \!\big(T+|{\cal I}m z|\big)^{\max(0,\ \!{\cal R}e(\omega-z)-2r)+\epsilon/2}\cr
\Big({1\over2}-a\le\ \!&{\cal R}e(z-\omega)+r\le 1-a+\log^{-1}\!V\ \ ,\ \ |w+r-\omega+z-1|\ge\log^{-1}\!T\Big)\ .}\leqno(16)$$

The lemma now follows from (15), (16) when $\ \!{\cal R}e(z-\omega)+r\le 1-a+\log^{-1}\!V\ ,\ \ \!|w+r\ \!+$

$-\ \!\omega+z-1|\ge\log^{-1}\!T\ \!\!,\ \!$ while it is trivial if $\ \!{\cal R}e(z-\omega)+r\ge 1-a+\log^{-1}\!V\ \!.$

Finally, it also holds for $\ |w+r-\omega+z-1|\le\log^{-1}\!T\ $ since $F_V(s)$ is olomorphic .

\qed

According to (10), put
$$\eqalign{&\quad\ \ F_V(s)\ =\ \sum_{n=1}^{+\infty}{c_n\over n^s}\qquad\ (\sigma>1)\cr
&c_n\ =\ c_n(V)\ =\sum_{d|n,\ d\le V}\mu(d)\ \!d^{1-s_V}\ .}\leqno(17)$$

Then, by (14)
$$G_{UV}(s)\ =\ {U\over2}^{s-s_0}\ \!\sum_{n=1}^{+\infty}\Big({c_n\over n^{s+iv}}+{c_n\over n^{s-iv}}\Big)\qquad\ (\sigma>1)\ .\leqno(18)$$

Also, the elementary upper bound
$$d(n)\ \!=\sum_{d|n}1\ \ \ll_{\epsilon}\ \!n^{\epsilon}\ \!\le V^{\epsilon}\qquad\big(n\le V\big)$$

gives 
$$V^{a-1}\sum_{d|n}\log d\ \le\ V^{a-1}\ \!d(n)\log n\ \le\ \!C(\epsilon)\ \!V^{a-1+2\epsilon}\ \!\le\ 1$$

if $\ \!0<\epsilon\le (1-a)/3\ \!$ and $\ \!V\ge V_0(\epsilon)\ \!.\ \!$ Hence, by lemma 2
$$c_n\ \!=\sum_{d|n}\mu(d)\ \!\Big(\exp\big\{(1-s_V)\log d\big\}-1\Big)\ \ll\ V^{a-1}\sum_{d|n}\log d\qquad\ \  \big(1<n\le V\big)\ \ .\leqno(19)$$

The crucial step in our method is the following result

{\bf Lemma 5.}\ \ {\it Let $\ \epsilon\ ,\ r\ ,\ a\ ,\ s_0\ ,\ T\ $ be as in (4) and let$\ G_{UV}(s)\ $be defined by (14). If} 
$$\eqalign{&\ \ \!2r\le\ \!z_0\ \!\le\ \!\min\big(3r\ \!,\ \!(2a-1)/10\ \!,\ \!(1-a)/5\big)\quad,\quad z_0\le z_1\le 2z_0\cr
&V\ge\max\big(T^{2/r},\ \!U\big)\quad,\quad U\ge\exp\big\{10/z_0\big\}\quad ,\quad J = 2 [z_0\log U+2]}$$

\qquad\qquad\quad\ \ {\it then}

$$\eqalign{&\ \ \!{(-z_0)^J\over J!}\ \!D^J G_{UV}(s_0+z_0-z_1)\ -\ {U^{z_0-z_1}(-z_0\log U)^J\over J!}\ \ll_{\epsilon}\ \!{U^{z_1}\over V^{2z_0-{r\over4}-6\epsilon}}\ \!+\ \!{U^{5z_0-r+{\epsilon\over2}}\over V^{5z_0}}\ +\cr
&\qquad\quad\ + \Big((VT)^{z_1-2r+4\epsilon} +\ \!V^{z_1-z_0-{7\over4}r+6\epsilon}+\ \!V^{z_1+z_0-{7\over4}r+6\epsilon}\ \!U^{-(1+2\log2)z_0}\Big)\ \!U^{-z_1}\ ,}$$
$$\eqalign{&\quad\!\sum_{1\le j<J}{(-z_0)^j\over j!}\ \!D^j G_{UV}(s_0)\ \!-\ \!\sum_{1\le j<J}{(-z_0\log U)^j\over j!}\ \!\ll_{\epsilon}\ \!V^{-2z_0+r}\ \!U^{z_0}\ \!+\ \!{U^{5z_0-r+\epsilon}\over V^{5z_0}}\ \!+\cr
&\qquad\quad\ \ +\Big((VT)^{z_0-2r+4\epsilon} +\ \!V^{2z_0-{7\over4}r+6\epsilon}\ \!U^{-(1+2\log2)z_0}\Big)U^{-z_0}\ \!+\ \!V^{{r\over4}+6\epsilon}\ \!U^{-r}\ .}$$

{\bf Proof.}\ \ Let
$$\delta(n)\ =\ \cases{ V^{a-1}\qquad\ \!if\quad\ 2\le n\le V\cr
 1\qquad\qquad if\qquad\ n> V\ \ \!.}\leqno(20)$$

Then, using (17), (19), (20) \big(see [3] Lemma 3.19\ with $\ a_n = c_n\ \!,\ x=Y\ \!,\ T=X\ \!,\ \alpha = 2\ \!,\ c =$

$=1-a+{\cal R} e\omega-r +\log^{-1}Y\ \!,\ \psi(n) =\delta(n)\ \!n^{\epsilon/2}$\big)
$$\eqalign{&\sum_{n\le Y}{c_n\over n^{s_0\pm iv-\omega}}\ =\ 1\ \!+\ \!{1\over2\pi i}\int_{1-a+{\cal R}e\omega-r +\log^{-1}\!Y -iX}^{1-a+{\cal R}e\omega-r +\log^{-1}\!Y+iX}\!\Big(F_V(s_0\pm iv-\omega+z)-1\Big)\ Y^z\ {dz\over z}\ \ \!+\cr
&\qquad\qquad\qquad\qquad\qquad\quad\ \ \!+\ O_{\epsilon}\bigg(Y^{{\cal R}e\omega-r-a+\epsilon}\Big({Y\over X} + \delta(2Y)\Big)\bigg)\ ,\cr
&\qquad\qquad\qquad\qquad\quad\ |\omega|\le z_1\quad,\quad 2\ \!\le\ \!Y\ \!\le\ \! X\quad,\quad X\ \!\ge\ \!\max(V,3T)\cr
&\qquad\ \ so\ \ that\qquad Y^{{\cal R}e\omega-r-a+\epsilon}\Big({Y\over X} + \delta(2Y)\Big)\ \!\ll\ \!V^{-a+{\cal R}e\omega-r+\epsilon}+ V^{a-1}\ \!\ll\ V^{a-1}}\leqno(21)$$

since, by hypothesis $\ \min\big(1-a\ \!,\ \!(2a-1)/2\big)\ge 5z_0\ge 5z_1/2\ge 2|{\cal R} e\omega|+r\ .$

Taking $\ \omega = z_1\ \!$ in (21) and appealing to lemma 3 with $\ \!s = s_0\pm iv-z_1+z \ \big(1-r\le{\cal R}e z\le$

$\le\ \!1-a +z_1-r+1/(\log Y)\ \!,\ |{\cal J}m z|\le  X\big)\ \!,\ \!$ we obtain  
$$\eqalign{&\quad\ {1\over2\pi i}\int_{1-a+z_1-r +\log^{-1}\!Y -iX}^{1-a+z_1-r +\log^{-1}\!Y+iX}\!\Big(F_V(s_0\pm iv-z_1+z)-1\Big)\ Y^z\ {dz\over z}\ \!=\ \!{1\over2\pi i}\bigg(\int_{z_1-r-iX}^{z_1-r+iX}\ +\cr
&+\ \!\int_{z_1-r+iX}^{1-a+z_1-r +\log^{-1}\!Y+iX}-\int_{z_1-r-iX}^{1-a+z_1-r +\log^{-1}\!Y -iX}\bigg)\Big(F_V(s_0\pm iv-z_1+z)-1\Big)\ Y^z\ {dz\over z}\ ,\cr
&\qquad\qquad\ \ {1\over2\pi i}\int_{z_1-r\pm iX}^{1-a+z_1-r +\log^{-1}\!Y\pm iX}\!\Big(F_V(s_0\pm iv-z_1+z)-1\Big)\ Y^z\ {dz\over z}\ \ll_{\epsilon}\cr
&\qquad\qquad\qquad\! \ll_{\epsilon}\ \!\Big(V^{-r+\epsilon}\ \!Y^{z_1-r} +\ \!V^{a-1+\epsilon}\ \!Y^{1-a+z_1-r}\Big)\ \!X^{2\epsilon-1}\ \!\ll\ V^{a-1}}\leqno(22)$$

and (21), (22) imply
$$\eqalign{&\sum_{n\le Y}{c_n\over n^{s_0\pm iv-z_1}}\ \!=\ \!1+\ \!{1\over2\pi i}\!\int_{z_1-r-iX}^{z_1-r+iX}\!\!\!\Big(F_V(s_0\pm iv-z_1+z)-1\Big)\ \!Y^z\ \!{dz\over z}\ \!+\ \!O\Big(V^{a-1}\Big)\cr
&\qquad\qquad\qquad\qquad\quad\ \big(2\le Y\le X\ \ ,\ \ X\ge\max(V,3T)\big)\ .}\leqno(23)$$

Now put $\ s_1 =\ \!s_0+z_0-z_1\ \!.\ $ From (14), (18), (23) we deduce
$$\eqalign{&{U\over 2}^{-z_1}\sum_{n\le Y}\Big({c_n\over n^{s_1+iv-z_0}}+{c_n\over n^{s_1-iv-z_0}}\Big)\ \!=\ \!{U\over 2}^{-z_1}\sum_{n\le Y}\Big({c_n\over n^{s_0+iv-z_1}}+{c_n\over n^{s_0-iv-z_1}}\Big)\ \!=\ U^{-z_1}\ \!+\cr
&\quad\ \ \!+\ \!{1\over2\pi i}\!\int_{z_1-r-iX}^{z_1-r+iX}\!\!\!\Big(G_{UV}(s_0-z_1+z)-U^{z-z_1}\Big)\ \!\Big({Y\over U}\Big)^z\ \!{dz\over z}\ +\ O\bigg(V^{a-1}\ \!U^{-z_1}\bigg)\ =\cr
&\ =\ \!U^{-z_1} +\ \!{1\over2\pi i}\!\int_{z_0-r-iX}^{z_0-r+iX}\!\!\!\Big(G_{UV}(s_0-z_1+z)-U^{z-z_1}\Big)\ \!\Big({Y\over U}\Big)^z\ \!{dz\over z}\ +\ O\bigg(V^{a-1}\ \!U^{-z_1}\bigg)}\leqno(24)$$

since, by lemma 4 \big(with $\omega =z_1$\big))

$$\eqalign{&\qquad\qquad{1\over2\pi i}\!\int_{z_0-r\pm iX}^{z_1-r\pm iX}\!\!\!\Big(G_{UV}(s_0-z_1+z)\ \!-U^{z-z_1}\Big)\ \!\Big({Y\over U}\Big)^z\ \!{dz\over z}\ \ll_{\epsilon}\cr
&\ll_{\epsilon}\ \!{(VX)^{\max(0,\ \!z_1-z_0-r)+4\epsilon}\over X U^{z_1}}\ \!Y^{z_1-r}\ \!\ll\ V^{-1+4z_0-3r+8\epsilon}\ \!U^{-z_1}\ \!\ll\ V^{a-1}\ \!U^{-z_1}\ .}$$

Moreover, using (3) and partial summation
$$\eqalign{&{1\over j!}\!\sum_{U<n\le X}{c_n\big(z_0\log(n/U)\big)^j\over n^{s_1\pm iv}}\ =\ {1\over U^{z_0}}\!\!\sum_{U<n\le X}{c_n\ \!p_j\big(z_0\log(n/U)\big)
\over n^{s_1\pm iv-z_0}}\ =\ {1\over U^{z_0}}\bigg(\sum_{U<n\le X}{c_n\over n^{s_1\pm iv-z_0}}\ \cdot\cr
&\qquad\qquad\quad\ \!\cdot p_j\big(z_0\log(X/U)\big)\ -\ z_0\!\int_U^X\!\!\sum_{U<n\le Y}{c_n\over n^{s_1\pm iv-z_0}}\ p'_j\big(z_0\log(Y/U)\big)\ {dY\over Y}\bigg)}\leqno(25)$$

for $\ j\ge1\ .\ $  Hence,  by (24), (25)
$$\eqalign{&\ {U\over2}^{z_0-z_1}\!\!\!\sum_{U<n\le X}\Big({c_n\over n^{s_1+iv}}+{c_n\over n^{s_1-iv}}\Big){\big(z_0\log(n/U)\big)^j\over j!}\ \!=\ {1\over2\pi i}\int_{z_0-r-iX}^{z_0-r+iX}\!\Big(G_{UV}(s_0-z_1 +z)\ +\cr
&\qquad -U^{z-z_1}\Big)\ \!\bigg(p_j\big(z_0\log(X/U)\big)\Big({X\over U}\Big)^z-\ z_0\!\int_U^Xp'_j\big(z_0\log(Y/U)\big)\Big({Y\over U}\Big)^z\ {dY\over Y}\bigg)\ \!{dz\over z}\ +\cr
&\qquad\qquad\ +\ \!{1\over U^{z_1}}\!\bigg(f(X)\ \!p_j\big(z_0\log(X/U)\big)\ \!-\ z_0\!\int_U^Xf(Y)\ \!p'_j\big(z_0\log(Y/U)\big)\ {dY\over Y}\bigg)\cr
&\quad\ where\qquad\quad\ \!f(u)\ \!\ll\ \!V^{a-1}\quad\!,\quad    1\le j\le J\quad ,\quad\!2\le U\le V\quad ,\quad X \ge\max(2V,3T)\ \ \!.}\leqno(26)$$

From now on we choose $\ \!X=V^2\ \!.\ \!$ According to (3), we then have, on integrating several times  
$$\eqalign{&\qquad\qquad\qquad\qquad\qquad\quad0 < p_j\big(z_0\log(V^2/U)\big) < 1\qquad\ \ \big(j\ge 0\big)\cr
&0<\!\sum_{j=0}^{J-1}p_j\big((z_0-{\cal R}ez)\log(V^2/U)\big)\ \!< \sum_{j=0}^Jp_j\big((z_0-{\cal R}ez)\log(V^2/U)\big)\ \!<\ \!1\qquad\ \big({\cal R}ez\le z_0-r\big)\cr
&\quad\!z_0\!\int_U^{V^2}\!p'_J\big(z_0\log(Y/U)\big)\ \!\Big({Y\over U}\Big)^z\ {dY\over Y}\ \!=\ \!z_0\!\int_1^{V^2/U}\!\!\Big(p_{J-1}(z_0\log y) - p_J(z_0\log y)\Big)\ \!y^z\ {dy\over y}\ =\cr
&\qquad\qquad\quad\ =\Big({z_0\over z_0-z}\Big)^{J+1}\bigg(\sum_{j=0}^J{\big((z_0-z)\log(V^2/U)\big)^j\over j!}\ \!\Big({V^2\over U}\Big)^{z-z_0}-\ 1\bigg)\ +\cr
&\qquad\qquad\quad\quad-\Big({z_0\over z_0-z}\Big)^J\bigg(\sum_{j=0}^{J-1}{\big((z_0-z)\log(V^2/U)\big)^j\over j!}\ \!\Big({V^2\over U}\Big)^{z-z_0}-\ 1\bigg)\ ,\cr
&\qquad\quad\! z_0\!\sum_{1\le j< J}\int_U^{V^2}\!p'_j\big(z_0\log(Y/U)\big)\ \!\Big({Y\over U}\Big)^z\ {dY\over Y}\ \!=\ \!z_0\!\sum_{1\le j<J}\int_1^{V^2/U}\!\!\Big(p_{j-1}(z_0\log y)\ +\cr
&\quad\quad\!- p_j(z_0\log y)\Big)\ \!y^z\ {dy\over y}\ \!=\ \!z_0\!\int_1^{V^2/U}\!\!\!\big(y^{-z_0}- p_{J-1}(z_0\log y)\big)\ \!y^z\ {dy\over y}\ =\ \!{z_0\over z_0-z}\ \!\Big(1\ \!+\cr
&\qquad\ -\!\big(V^2/U\big)^{z-z_0}\Big)\ \! +\ \!\Big({z_0\over z_0-z}\Big)^J\bigg(\sum_{j=0}^{J-1}{\big((z_0-z)\log(V^2/U)\big)^j\over j!}\ \!\Big({V^2\over U}\Big)^{z-z_0}-\ 1\bigg)\ .}\leqno(27)$$

On the other hand 
$$\eqalign{&\quad\ \! {1\over2\pi i}\int_{z_0-r-iV^2}^{z_0-r+iV^2}\!\Big(G_{UV}(s_0-z_1+z)-U^{z-z_1}\Big)\ \!\big(V^2/U\big)^z\ \!{dz\over z}\ =\ G_{UV}(s_0-z_1)-U^{-z_1}\ +\cr
&+\ \!{1\over2\pi i}\Bigg(\int_{-z_0-iV^2}^{-z_0+iV^2}\!+\!\int_{-z_0+iV^2}^{z_0-r+iV^2}-\int_{-z_0-iV^2}^{z_0-r-iV^2}\Bigg)\ \!\Big(G_{UV}(s_0-z_1+z)-U^{z-z_1}\Big)\ \!\big(V^2/U\big)^z\ \!{dz\over z}\ .}\leqno(28)$$

Using (10), (14), (15) and lemma 4 (with $\ \!\omega = z_1),\ \!$ we now obtain
$$\eqalign{&\ G_{UV}(s_0-z_1)-U^{-z_1}\ \!\ll_{\epsilon}(VT)^{z_1-2r+4\epsilon}\ \!U^{-z_1}\quad,\quad \int_{-z_0\pm iV^2}^{z_0-r\pm iV^2}\Big(G_{UV}(s_0-z_1+z)-U^{z-z_1}\Big)\ \!\cdot\cr
&\qquad\quad\cdot\big(V^2/U\big)^z\ \!{dz\over z}\ \!\ll_{\epsilon}\ \!\big(V^{3z_1+z_0-6r+12\epsilon} +\ \!V^{\max(2z_0-2r,\ \!3z_1-z_0-5r)+12\epsilon}\big)\ \!U^{-z_1} V^{-2}\ \!,\cr
&\int_{-z_0-iV^2}^{-z_0+iV^2}\!\ll_{\epsilon} {V^{\max(0,\ \!z_1+z_0-2r)-2z_0+5\epsilon}\over U^{z_1}}\bigg(\int_{-V^2}^{V^2}\Big(1+\big|\zeta\big(a-z_1-z_0+r+i(v+y)\big)\big|^2\Big)\ \! {dy\over|y|+r}\bigg)^{1\over2}}\leqno(29)$$

where, by hypothesis, $\ \ \!2r\le z_0\le z_1\le 2z_0\le\min\big(6r\ \!,\ \!(2a-1)/5\ \!,\ \!2(1-a)/5\big)\le 3/50\ .\ $ Hence
$$\eqalign{&\ z_1+z_0-2r>0\ \ ,\ \ V^{\max(3z_1+z_0-6r\ \!,\ \!2z_0-2r\ \!,\ \!3z_1-z_0-5r)+12\epsilon-2}\ \!\le\ V^{15r+12\epsilon-2}\ \!\le\ V^{-{37\over20}+12\epsilon}\ \!,\cr
&1-10r\ \!\ge\ \!a\ \!\ge\ \!a-z_1-z_0+r\ \!=\ \!{1\over2} +{2a-1\over 2} -z_1-z_0+r\ \!\ge\ \!{1\over2} +4z_0-z_1+r\ \!\ge\ \!{1\over2}+5r\ \!.}$$

Putting $\ L\ \!=\ \!\Big[{\log(V^2/T
)\over\log2}\Big]\ $ we then have \big(see [3] Theorem 7.2 (A) with $\ \!\sigma\ge {1\over2}+5r$\big)
$$\eqalign{&\qquad\qquad\qquad\quad\ \int_{-V^2}^{V^2}\big|\zeta\big(a-z_1-z_0+r+i(v+y)\big)\big|^2\ \!{dy\over|y|+r}\ \!\ll\ {1\over r}\ +\cr
&+\!\int_2^{2T}\!\big|\zeta(a-z_1-z_0+r+it)\big|^2\ \!dt\ \!+\ \!{1\over T} \sum_{\ell=1}^L\ \!{1\over 2^{\ell}}\!\int_{2^{\ell}T}^{2^{\ell+1}T}\big|\zeta(a-z_1-z_0+r+it)\big|^2\ \!dt\ \ll\cr
&\qquad\qquad\qquad\qquad\ \ll\ {1\over r} + {T + L\over r}\ \ll\ {T+\log V\over r}\ \ll_{\epsilon}\ T+\log V\ .}\leqno(30)$$

Inserting the above estimates in equation (28), we see that
$$\eqalign{&{1\over2\pi i}\int_{z_0-r-iV^2}^{z_0-r+iV^2}\!\Big(G_{UV}(s_0-z_1+z)-U^{z-z_1}\Big)\ \!\big(V^2/U\big)^z\ \!{dz\over z}\ \ll_{\epsilon}\cr
&\qquad\quad\ll_{\epsilon}\ \!\Big((VT)^{z_1-2r+4\epsilon}+\ \!V^{z_1-z_0-2r+6\epsilon}\ \!T^{1/2}\Big)\ \!U^{-z_1}} $$

and (27) implies
$$\eqalign{&\quad\ {1\over2\pi i}\int_{z_0-r-iV^2}^{z_0-r+iV^2}\!\Big(G_{UV}(s_0-z_1+z)-U^{z-z_1}\Big)\ \!p_J\big(z_0\log(V^2/U)\big)\ \!\big(V^2/U\big)^z\ \!{dz\over z}\ \ll_{\epsilon}\cr
&\qquad\qquad\qquad\quad\ \ \ll_{\epsilon}\ \!\Big((VT)^{z_1-2r+4\epsilon}+\ \!V^{z_1-z_0-2r+6\epsilon}\ \!T^{1/2}\Big)\ \!U^{-z_1}\ \!,\cr
&{1\over2\pi i}\int_{z_0-r-iV^2}^{z_0-r+iV^2}\!\Big(G_{UV}(s_0-z_0+z)-U^{z-z_0}\Big)\!\sum_{1\le j<J}p_j\big(z_0\log(V^2/U)\big)\ \!\big(V^2/U\big)^z\ \!{dz\over z}\ \ll_{\epsilon}\cr
&\qquad\qquad\qquad\qquad\ \ \ll_{\epsilon}\ \!\Big((VT)^{z_0-2r+4\epsilon}\ \!+\ \!V^{-2r+6\epsilon}\ \!T^{1/2}\Big)\ \!U^{-z_0}\ .}\leqno(31)$$

Here \big(and in (35), (37), (38), (39)\big) we take $\ \!z_1=z_0\ \!$ in the latter formula.

Also, by the same argument as before 
$$\eqalign{&\ {1\over2\pi i}\int_{z_0-r-iV^2}^{z_0-r+iV^2}\!\Big(G_{UV}(s_0-z_1+z)-U^{z-z_1}\Big)\ \!\Big({z_0\over z_0-z}\Big)^{J+1}\sum_{j=0}^J{\big((z_0-z)\log(V^2/U)\big)^j\over j!}\ \cdot\cr
&\ \cdot\Big({V^2\over U}\Big)^{z-z_0}\ {dz\over z}\ \!=\ \!\Big(G_{UV}(s_0-z_1)-U^{-z_1}\Big)\!\sum_{j=0}^Jp_J\big(z_0\log(V^2/U)\big)\ \!+\ \!{1\over2\pi i}\Bigg(\!\int_{-z_0-iV^2}^{-z_0+iV^2}\!\! +\cr
&\qquad\qquad+\!\int_{-z_0+iV^2}^{z_0-r+iV^2}\!-\!\int_{-z_0-iV^2}^{z_0-r-iV^2}\Bigg) \Big(G_{UV}(s_0-z_1+z)-U^{z-z_1}\Big)\Big({z_0\over z_0-z}\Big)^{J+1}\cdot\cr
&\quad\cdot\ \!\!\sum_{j=0}^J{\big((z_0-z)\log(V^2/U)\big)^j\over j!}\Big({V^2\over U}\Big)^{z-z_0}\ \! {dz\over z}\ \ll_{\epsilon}\ \!{\ \!(VT)^{z_1-2r+4\epsilon}\over U^{z_1}}\ \!+\ \!{\big(V/U\big)^{z_1+z_0}\ \!T^{1/2}\over 2^J\ \!V^{2r-6\epsilon}}}\leqno(32)$$

the above upper bound being a consequence of the following inequalities \big(see (3), (27)\big)
$$\eqalign{&\qquad\qquad\qquad\qquad\qquad\bigg|{\big(V^2/U\big)^{z-z_0}\over(z_0-z)^{J+1}}\sum_{j=0}^J{\big((z_0-z)\log(V^2/U)\big)^j\over j!}\bigg|\ \le\cr
&\le\ {1\over(z_0-{\cal R}ez)^{J+1}}\sum_{j=0}^Jp_j\big((z_0-{\cal R}ez)\log(V^2/U)\big)\ \le\ {1\over(z_0-{\cal R}ez)^{J+1}}\qquad\ \big({\cal R}ez\le z_0-r\big)}$$

and, appealing to (15) and to lemma 4 \big(see (27), (29), (30), (31)\big), of these below
$$\eqalign{&\qquad\qquad\quad\Big(G_{UV}(s_0-z_1)-U^{-z_1}\Big)\!\sum_{j=0}^Jp_J\big(z_0\log(V^2/U)\big)\ \!\ll_{\epsilon}\ \!{(VT)^{z_1-2r+4\epsilon}\over U^{z_1}}\ \!,\cr
&\qquad\int_{-z_0-iV^2}^{-z_0+iV^2}\big|G_{UV}(s_0-z_1+z)-U^{z-z_1}\Big|\ \Big|{z_0\over z_0-{\cal R}e z}\Big|^{J+1}\ \!\Big|{dz\over z}\Big|\ \ll_{\epsilon}\ \!{V^{z_1+z_0-2r+5\epsilon}\over 2^J\ \!U^{z_1+z_0}}\ \!\cdot\cr
&\quad\cdot\bigg(\int_{-V^2}^{V^2}\Big(1+\big|\zeta\big(a-z_1-z_0+r+i(v+y)\big)\big|^2\Big)\ \! {dy\over|y|+r}\bigg)^{1\over2}\ \!\ll_{\epsilon}\ \!{V^{z_1+z_0-2r+6\epsilon}\ \!T^{1/2}\over 2^J\ \!U^{z_1+z_0}}\ \ \!,\cr
&\int_{-z_0\pm iV^2}^{z_0-r\pm iV^2}\big|G_{UV}(s_0-z_1+z)-U^{z-z_1}\Big|\ \Big|{z_0\over z_0-{\cal R}e z}\Big|^{J+1}\ \!\Big|{dz\over z}\Big|\ \ll_{\epsilon}\ \!{V^{3(z_1+z_0-2r)+12\epsilon}{z_0}^J\over V^2\ \!U^{z_1-z_0+r}\ \!{r}^J}\ \!\ll\cr
&\ll\ {V^{3(z_1+z_0-2r)+12\epsilon-2}\over 2^J\ \!U^{z_1-(1+2\log6)z_0+r}}\ =\ {V^{z_1+z_0-2r+6\epsilon}\over 2^J\ \!U^{z_1+z_0}}\cdot {U^{2(1+\log6)z_0-r}\over V^{2(1-z_1-z_0+2r-3\epsilon)}}\ <\ {V^{z_1+z_0-2r+6\epsilon}\ \!T^{1/2}\over 2^J\ \!U^{z_1+z_0}}}$$

since, by hypothesis, $\ \!\big(z_0/r\big)^J\!\le 3^J\le 6^{2z_0\log U+4}\ \!2^{-J}\ll\ \!U^{2z_0\log6}\ \!2^{-J}\ \!$ and $\ U^{2(1+\log6)z_0-r}\le$

$\le\ \!V^{2(1+\log6)z_0-r}\le\ \!V^{2(1-z_1-z_0+2r-3\epsilon)}\ .$

Quite similarly we obtain
$$\eqalign{&\ \ {1\over2\pi i}\int_{z_0-r-iV^2}^{z_0-r+iV^2}\!\!\Big(G_{UV}(s_0-z_1+z)-U^{z-z_1}\Big)\ \!\Big({z_0\over z_0-z}\Big)^J\sum_{j=0}^{J-1}{\big((z_0-z)\log(V^2/U)\big)^j\over j!}\ \cdot\cr
&\qquad\ \ \cdot\Big({V^2\over U}\Big)^{z-z_0}\ \!{dz\over z}\ \ll_{\epsilon}\ \!(VT)^{z_1-2r+4\epsilon}\ \!U^{-z_1} +\ \!\big(V/U\big)^{z_1+z_0}\ \!V^{-2r+6\epsilon}\ \!2^{-J}\ \!T^{1/2}\ \!,\cr
&{1\over2\pi i}\int_{z_0-r-iV^2}^{z_0-r+iV^2}\!\Big(G_{UV}(s_0-z_1+z)-U^{z-z_1}\Big)\Big({z_0\over z_0-z}\Big)^{\ell}\ \!{dz\over z}\ \!=\ \!G_{UV}(s_0-z_1)-U^{-z_1}\ +\cr
&\qquad\qquad\quad+\ \!O\Big({\big(V/U\big)^{z_1+z_0}T^{1/2}\over2^J V^{2r-6\epsilon}}\Big)\ \!\ll\ \!{(VT)^{z_1-2r+4\epsilon}\over U^{z_1}}+{\big(V/U\big)^{z_1+z_0}T^{1/2}\over2^J\ \!V^{2r-6\epsilon}}}\leqno(33)$$

for $\ \!\ell = J,\ \!J+1 .\ \!$

Furthermore (15) (with $\ \!\omega = z_1)\ \!$ implies \big(see (28), (29)\big)
$$\eqalign{&\quad\ \ {1\over2\pi i}\int_{z_0-r-iV^2}^{z_0-r+iV^2}\!\Big(G_{UV}(s_0-z_1+z)-U^{z-z_1}\Big)\ \!{z_0\over z_0-z}\ \!\Big(1-\ \!\big(V^2/U\big)^{z-z_0}\Big)\ \! {dz\over z}\ \ll_{\epsilon}\cr
&\ll_{\epsilon}\ \!V^{\max(0,\ \!z_1-z_0-r)+5\epsilon}\ \!U^{z_0-z_1-r}\bigg(\int_{-V^2}^{V^2}\Big(1+\big|\zeta\big(a-z_1+z_0+i(v+y)\big)\big|^2\Big)\ \! {dy\over|y|+r}\bigg)^{1\over2}}$$

and then, observing that $1-10r\ge a\ge a-z_1+z_0 = 1/2 + (2a-1)/2-z_1+z_0\ge 1/2-z_1\ \!+\ \!$

$+\ \!6z_0\ge 1/2 +8r\ \!,\ $ we get \big(see (29), (30), (31)\big)
$$\eqalign{&{1\over2\pi i}\int_{z_0-r-iV^2}^{z_0-r+iV^2}\!\Big(G_{UV}(s_0-z_1+z)-U^{z-z_1}\Big)\ \!{z_0\over z_0-z}\ \!\Big(1-\ \!\big(V^2/U\big)^{z-z_0}\Big)\ \! {dz\over z}\ \ll_{\epsilon}\cr
&\qquad\qquad\qquad\qquad\ll_{\epsilon}\ \!\Big(1+V^{z_1-z_0-r}\Big)\ \!V^{6\epsilon}\ \!U^{z_0-z_1-r}\ \!T^{1/2}\ .}\leqno(34)$$

From (27), (32), (33), (34) it follows that
$$\eqalign{&\quad{z_0\over2\pi i}\int_{z_0-r-i V^2}^{z_0-r+iV^2}\!\Big(G_{UV}(s_0-z_1 +z)-U^{z-z_1}\Big)\ \!{dz\over z}\!\int_U^{V^2}\!\!p'_J\big(z_0\log(Y/U)\big)\ \!\Big({Y\over U}\Big)^z\ \!{dY\over Y}\ \ll_{\epsilon}\cr
&\qquad\qquad\qquad\qquad\qquad\ \ll_{\epsilon}\ {(VT)^{z_1-2r+4\epsilon}\over U^{z_1}}\ \!+\ \!{V^{z_1+z_0-{7\over4}r+6\epsilon}\over U^{z_1+(1+2\log2)z_0}}\ ,\cr
&\ {z_0\over2\pi i}\sum_{1\le j<J} \int_{z_0-r-iV^2}^{z_0-r+iV^2}\!\Big(G_{UV}(s_0-z_0+z)-U^{z-z_0}\Big)\ \!{dz\over z}\int_U^{V^2}\!p'_j\big(z_0\log(Y/U)\big)\Big({Y\over U}\Big)^z\ {dY\over Y}\ \!\ll_{\epsilon}\cr
&\qquad\qquad\ \ \ll_{\epsilon}\ \!(VT)^{z_0-2r+4\epsilon}\ \!U^{-z_0}\ \!+\ V^{2z_0-{7\over4}r+6\epsilon}\ \!U^{-2(1+\log2)z_0}\ \!+\ V^{{r\over 4}+6\epsilon}\ \!U^{-r}}\leqno(35)$$

since $\ \!J>2z_0\log U\ \!$ and $\ \!T\le V^{r/2}\ \!.\ $ Here we point out that, according to (27), the upper bound

(34) appears only in the latter estimate, when we sum over $j$ and $\ \!z_1=z_0\ .$

Finally, by (3), (26), (27)
$$\eqalign{&\quad{f(V^2)\over U^{z_1}}\ \!p_J\big(z_0\log(V^2/U)\big)\ \!\ll\ \!V^{a-1}\ \!U^{-z_1}\quad\!,\quad{f(V^2)\over U^{z_0}}\sum_{1\le j<J}p_j\big(z_0\log(V^2/U)\big)\ \!\ll\ \!V^{a-1}\ \!U^{-z_0}\ \!,\cr
&{z_0\over U^{z_1}}\!\int_U^{V^2} p'_J\big(z_0\log(Y/U)\big)\ \!f(Y)\ {dY\over Y}\ =\ \!{1\over U^{z_1}}\!\int_0^{z_0\log(V^2/U)}\!f(Ue^{u/z_0})\ \Big({u^{J-1}\over(J-1)!}-{u^J\over J!}\Big)\ \! e^{-u}\ \!du\ \!\ll\cr
&\qquad\qquad\qquad\qquad\!\ll\ V^{a-1}\ \!U^{-z_1}\!\!\int_0^{+\infty}\!\!\Big({u^{J-1}\over(J-1)!}+{u^J\over J!}\Big)\ \!{du\over e^u}\ \ll\ V^{a-1}\ \!U^{-z_1}\ \!,\cr
&\ \ \!{z_0\over U^{z_0}}\!\int_U^{V^2}\!\!\sum_{1\le j<J} p'_j\big(z_0\log(Y/U)\big)\ \!f(Y)\ {dY\over Y}\ \ll\ {V^{a-1}\over U^{z_0}}\!\int_0^{+\infty}\!\!\Big(1+{u^{J-1}\over(J-1)!}\Big)\ \!{du\over e^u}\ \ll\ V^{a-1}\ \!U^{-z_0}\ \! .  }\leqno(36)$$

Inserting the estimates (31), (35), (36) in the basic identity (26), we now obtain
$$\eqalign{&\ {U\over2}^{z_0-z_1}\!\!\sum_{U<n\le V^2}\Big({c_n\over n^{s_1+iv}}+{c_n\over n^{s_1-iv}}\Big){\big(z_0\log(n/U)\big)^J\over J!}\ \ll_{\epsilon}\ \!(VT)^{z_1-2r+4\epsilon}\ \!U^{-z_1}\ \!+\cr
&\qquad\qquad\! +\ V^{z_1-z_0-{7\over4}r+6\epsilon}\ \!U^{-z_1}\ \!+\ V^{z_1+z_0-{7\over4}r+6\epsilon}\ \!U^{-z_1-(1+2\log2)z_0}\ ,\cr
&{1\over2}\!\sum_{1\le j<J}\ \!\sum_{U<n\le V^2}\Big({c_n\over n^{s_0+iv}}+{c_n\over n^{s_0-iv}}\Big){\big(z_0\log(n/U)\big)^j\over j!}\ \ll_{\epsilon}\ \!(VT)^{z_0-2r+4\epsilon}
\ \!U^{-z_0}\ \!+\cr
&\qquad\qquad\qquad\quad\ +\ V^{2z_0-{7\over4}r+6\epsilon}\ \!U^{-2(1+\log2)z_0}\ \!+\ V^{{r\over 4}+6\epsilon}\ \!U^{-r}}\leqno(37)$$

since, by hypothesis, $\ \!2r\le z_0\le z_1\le 2z_0\le 2(1-a)/5\ \!$ and where $\ \!s_1 = s_0+z_0-z_1\ $\big(see (24)\big) . 

Also
$$\eqalign{&\qquad\qquad\qquad\qquad\quad 1-s_0=1-a-r\ge 5z_0-r\cr
&s_0+z_1\le a+r+{2(1-a)\over5} = 1+r-{3(1-a)\over5}\le 1-3z_0+r\le 1-5r\ .}$$

Therefore, when $\ \!|z|\le z_1\ \!,\ \!$ (19) yields
$${1\over U^z}\!\!\!\sum_{1<n\le U}{c_n\over n^{s_0\pm iv-z}}\ \ll\ {U^{z_1}\over V^{1-a}}\sum_{d\le U}\ \!{\log d\over d^{s_0+z_1}}\!\sum_{\ell\le U/d}{1\over\ell^{s_0+z_1}}\ \ll\ {U^{1-s_0}\over r\ \!V^{1-a}}\!\sum_{d\le U}{\log d\over d}\ \ll_{\epsilon}\ \!{U^{5z_0-r+{\epsilon\over2}}\over V^{5z_0}}\ .$$

We now apply Cauchy's inequality for the coefficients of a power series to the function
$$\Phi_1(z)\ \!=\ \! {1\over2}\!\sum_{1\le n\le U}\Big({c_n\over n^{s_0+iv}}+{c_n\over n^{s_0-iv}}\Big)\Big({U\over n}\Big)^{-z}$$

in the circle $\ \!|z-z_1+z_0| \le z_0\ \!$ (so that $\ \!|z|\le |z-z_1+z_0|\ \!+z_1-z_0\le z_1)\ \!.\ $ From the above

estimate we then obtain
$$\eqalign{&{U\over2}^{z_0-z_1}\!\!\!\sum_{1\le n\le U}\Big({c_n\over n^{s_1+iv}}+{c_n\over n^{s_1-iv}}\Big){\big(z_0\log(n/U)\big)^J\over J!}\ \!=\ {U^{z_0-z_1}(-z_0\log U)^J\over J!}\ \!+\ \!O_{\epsilon}\bigg({U^{5z_0-r+{\epsilon\over 2}}\over V^{5z_0}}\bigg)\cr
&\ {1\over2}\!\sum_{1\le j<J}\ \!\sum_{1\le n\le U}\Big({c_n\over n^{s_0+iv}}+{c_n\over n^{s_0-iv}}\Big){\big(z_0\log(n/U)\big)^j\over j!}\ \!=\!\sum_{1\le j<J}{(-z_0\log U)^j\over j!}\ \!+\ \!O_{\epsilon}\bigg({U^{5z_0-r+\epsilon}\over
V^{5z_0}}\bigg)}\leqno(38)$$

recalling that in the latter equation one takes $\ \!z_1=z_0\ \!$ and that $\ \! J\ll z_0\log U\ll_{\epsilon} U^{\epsilon/2}\ \!.\ \!$

On the other hand, when $\ |\omega-z_1+z_0|\le z_0\ \!,\ $ we deduce from (14), (21)
$$\eqalign{&{U\over 2}^{-\omega}\!\!\!\sum_{n\le V^2}\!\!\Big({c_n\over n^{s_0+iv-\omega}}+{c_n\over n^{s_0-iv-\omega}}\Big) =\ {1\over U^{\omega}} + {1\over2\pi i}\!\int_{1-a+{\cal R}e\omega-r +(2\log V)^{-1} -iV^2}^{1-a+{\cal R}e\omega-r +(2\log V)^{-1}\!+iV^2} \!\!\!\!\Big(G_{UV}(s_0-\omega+z)\ \!+\cr
&\ \!-U^{z-\omega}\Big)\ \!\Big({V^2\over U}\Big)^z\ \!{dz\over z}\ \!+\ \!O\bigg({U^{-{\cal R}e\omega}\over V^{1-a}}\bigg)\ \!=\ G_{UV}(s_0-\omega)\ \!+\ \!{1\over2\pi i}\bigg(\int_{-z_0+iV^2}^{1-a+{\cal R}e\omega-r+(2\log V)^{-1}+iV^2}\!+\cr
&-\!\!\int_{-z_0-iV^2}^{1-a+{\cal R}e\omega-r+(2\log V)^{-1}-iV^2}\!\!\!+\!\int_{-z_0-iV^2}^{-z_0+iV^2}\bigg)\Big(G_{UV}(s_0-\omega+z) -U^{z-\omega}\Big)\ \!{V^{2z}\over U^z}\ \!{dz\over z} + O\bigg({U^{z_1}\over V^{5z_0}}\bigg)\ }$$

since, as we observed, $\ \!1-a\ge 5z_0\ \!$ and $\ |{\cal R}e\omega|\le|\omega|\le|\omega+z_0-z_1|+z_1-z_0\le z_1\ \!.\ $

By (15) and lemma 4 we have further \big(see (28), (29), (30). (31)\big)
$$\eqalign{&\qquad\quad\ \int_{-z_0-iV^2}^{-z_0+iV^2}\Big(G_{UV}(s_0-\omega+z) -U^{z-\omega}\Big)\ \!{V^{2z}\over U^z}\ \!{dz\over z}\ \!\ll_{\epsilon}\ \!{V^{\max(0,\ \!{\cal R} e\omega+z_0-2r)-2z_0+5\epsilon}\over U^{{\cal R} e\omega}}\ \!\cdot\cr
&\cdot\bigg(\int_{-V^2}^{V^2}\Big(1+\big|\zeta\big(a-{\cal R} e\omega-z_0+r+i(v+y)\big)\big|^2\Big)\ \!{dy\over|y|+r}\bigg)^{1\over2} \ll_{\epsilon}\ \!V^{6\epsilon}\ \!T^{1\over2}\Big({U^{z_1}\over V^{2z_0}} + {V^{z_1-z_0-2r}\over U^{z_1}}\Big)\ ,\cr
&\ \ {1\over2\pi i}\int_{-z_0\pm iV^2}^{1-a+{\cal R}e\omega-r+(2\log V)^{-1}\pm iV^2}\!\!\!\Big(G_{UV}(s_0-\omega+z) -U^{z-\omega}\Big)\ \!\big(V^2/U\big)^z\ \!{dz\over z}\ \!\ll_{\epsilon} \big(U^{z_1}\ \!V^{12\epsilon-2z_0}\ +\cr
&\qquad\qquad\qquad\qquad +\ \!V^{3z_1+z_0-6r+12\epsilon}\ \!U^{-z_1}+\ \!V^{2(1-a)+2z_1-2r+12\epsilon}\ \!U^{-z_1}\big)V^{-2}\ \!.}$$

Then, using the inequality $\ T\le V^{r/2}$
$$\ \ {U\over 2}^{-\omega}\!\!\!\sum_{n\le V^2}\!\!\Big({c_n\over n^{s_0+iv-\omega}}+{c_n\over n^{s_0-iv-\omega}}\Big)\ \!-\ \!G_{UV}(s_0-\omega)\ \ll_{\epsilon}\ \Big({U^{z_1}\over V^{2z_0}}\ \!+\ \!{V^{z_1-z_0-2r}\over U^{z_1}}\Big)V^{{r\over4}+6\epsilon}$$

uniformly for $\ |\omega+z_0-z_1|\le z_0\ .$

Therefore, by applying Cauchy's inequality to the function
$$\Phi_2(\omega)\ \!=\ \! {1\over2}\!\sum_{1\le n\le V^2}\Big({c_n\over n^{s_0+iv}}+{c_n\over n^{s_0-iv}}\Big)\Big({U\over n}\Big)^{-\omega} -\  G_{UV}(s_0-\omega)$$

in the circle  $\ \!|\omega-z_1+z_0| \le z_0\ \!$ (as we did before)
$$\eqalign{&{U\over2}^{z_0-z_1}\!\!\sum_{1\le n\le V^2}\Big({c_n\over n^{s_1+iv}}+{c_n\over n^{s_1-iv}}\Big){\big(z_0\log(n/U)\big)^J\over J!}\ \!=\ \!{(-z_0)^J\over J!}\ \!D^J\ \! G_{UV}(s_0+z_0-z_1)\ +\cr
&\qquad\qquad\qquad\quad+\ O_{\epsilon}\Big(V^{-2z_0+{r\over4}+6\epsilon}\ \!U^{z_1} +\ \!V^{z_1-z_0-{7\over4}r+6\epsilon}\ \!U^{-z_1}\Big)\ ,\cr
&\ \ {1\over2}\!\sum_{1\le j<J}\ \!\sum_{1\le n\le V^2}\Big({c_n\over n^{s_0+iv}}+{c_n\over n^{s_0-iv}}\Big){\big(z_0\log(n/U)\big)^j\over j!}\ \!=\ \!\sum_{1\le j<J}{(-z_0)^j\over j!}\ \!D^j\ \! G_{UV}(s_0)\ +\cr
&\qquad\qquad\qquad\qquad\quad\ \ +\ O\Big( V^{-2z_0+r}\ \!U^{z_0} +\ \!V^{-r}\ \!U^{-z_0}\Big)\ .}$$

Once again we point out that in the latter equation we take $\ \!z_1=z_0\ \!$ and use the upper bound

$\ \! J\ll z_0\log U\ll V^{{3\over4}r-6\epsilon}\ \ \big(0< \epsilon\le r/100\ \!$ by (4)\big) .

Hence
$$\eqalign{&\quad D^J\ \! G_{UV}(s_0+z_0-z_1)\ \!=\ \!{U\over2}^{z_0-z_1}\!\!\sum_{1<n\le V^2}\Big({c_n\over n^{s_1+iv}}+{c_n\over n^{s_1-iv}}\Big){\big(z_0\log(n/U)\big)^J\over J!}\ +\cr
&\qquad+{U^{z_0-z_1}(-z_0\log U)^J\over J!}\ \!+\ O_{\epsilon}\Big(V^{-2z_0+{r\over4}+6\epsilon}\ \!U^{z_1} +\ \!V^{z_1-z_0-{7\over4}r+6\epsilon}\ \!U^{-z_1}\Big)\ ,\cr
&\sum_{1\le j<J}{(-z_0)^j\over j!}\ \!D^j\ \! G_{UV}(s_0)\ \!=\ \!{1\over2}\!\sum_{1\le j<J}\ \!\sum_{1<n\le V^2}\Big({c_n\over n^{s_0+iv}}+{c_n\over n^{s_0-iv}}\Big){\big(z_0\log(n/U)\big)^j\over j!}\ +\cr
&\qquad\qquad\quad\ \!+\ \!\sum_{1\le j<J}{(-z_0\log U)^j\over j!}\ \!+\  O\Big( V^{-2z_0+r}\ \!U^{z_0} +\ \!V^{-r}\ \!U^{-z_0}\Big)\ .}\leqno(39)$$

The lemma now follows from (37), (38), (39) .

\qed

\non

Let $\ b\ $ and $\ \rho_0 =\beta_0+i\gamma_0\ $ be as in (1) and (2), respectively. According to (1), (2), (4)

put
$$\eqalign{&\qquad\ \ \!T=2\gamma_0/3\quad,\quad r = \min\Big({1\over100}\ \!,\ {20b\over181}\ \!,\ {10(1-2b)\over221}\Big)\cr
&\epsilon\le r/100\quad,\quad a = b + (1+2r)/2\quad,\quad w = a + iv = a + i\gamma_0\cr
&\quad  z_0 = w+r-\rho_0 = s_0 - \beta_0\qquad\Big(b+{1\over2}-\epsilon\le \beta_0\le b+{1\over2}\Big)}\leqno(40)$$

whence \big(see (14)\big)
$$\eqalign{&\qquad\qquad\ 2r\le z_0\le 201r/100\le\min\big((2a-1)/10\ \!,\ \!(1-a)/5\big)\cr
&G_{UV}(s_0-z_0)\ \!=\ \!U^{-z_0}\ \!{\cal R}e\big\{F_V(w+r-z_0)\big\}\ \!=\ \!U^{-z_0}\ \!{\cal R}e\big\{F_V(\rho_0)\big\}\ \!=\ \!0}\leqno(41)$$
 
while lemma 3 implies, when $\ \! V\ge T^{2/r}$
$$G_{UV}(s_0)\ \!=\ \!{\cal R}e\big\{F_V(w+r)\big\} = 1 + O_{\epsilon}\Big(V^{\epsilon-2r}\ \!T^{\epsilon}\Big) = 1 +  O_{\epsilon}\Big(V^{-2(r-\epsilon)}\Big)\leqno(42)$$

Take further $\ \!U = V^{2/3}\ \!$ and recall that $\ \!G_{UV}(s)\ \!$ is real for real$\ s\ \!.\ \!$ Then,  by Taylor's formula,

given any integer $J>1\ \!,\ \!$ we can find a number $\ \!s_1\ \!,\ \ \!s_0-z_0<s_1< s_0\ \!,\ $ such that

$$G_{UV}(s_0-z_0)\ \!=\ \!G_{UV}(s_0)\ \!+\ \!\sum_{j=1}^{J\!-1}{(-z_0)^j\over j!}\ \!D^j\ \!G_{UV}(s_0)\ \!+\ \!{(-z_0)^J\over J!}\ \!D^J\ \!G_{UV}(s_1)\ .\leqno(43)$$

Choose $\ \!J = 2[(z_0\log U + 2]\ \!$ and apply lemma 5 with $\ \!z_1=z_0+s_0-s_1\ \!.\ $ It follows from (40),

(41), (42), (43)
$$\eqalign{&\ -1 + O_{\epsilon}\Big(V^{-2(r-\epsilon)}\Big)\ \!=\ \!\sum_{j=1}^{J-1}{(-z_0)^j\over j!}\ \!D^j\ \!G_{UV}(s_0)\ \!+\ \!{(-z_0)^J\over J!}\ \!D^J\ \!G_{UV}(s_1)\ \!=\sum_{j=1}^{J-1}{(-z_0\log U)^j\over j!}\ \!+\cr
&\quad\ \ +\ \!{U^{z_0-z_1}(-z_0\log U)^J\over J!}\ \!+\ \!O_{\epsilon}\bigg(V^{{2\over3}z_1-2z_0+{r\over4}+6\epsilon}+\ \!V^{-{5\over3}z_0-{2\over3}r+{1\over 3}\epsilon} +\ \!V^{{z_1\over3}-2r+4\epsilon}\ \!T^{z_1-2r+4\epsilon}\ \!+\cr
&+ V^{{z_1\over3}-z_0-{7\over4}r+6\epsilon} +\ \!V^{{z_1+(1-4\log2)z_0\over3}-{7r\over4}+6\epsilon} + V^{-{4\over3}z_0+r} + V^{-{5\over3}z_0-{2\over3}r+{2\over 3}\epsilon} + V^{{z_0\over3}-2r+4\epsilon}\ \!T^{z_0-2r+4\epsilon}\ \!+\cr
&+\Big(V^{{(2-4\log2)z_0\over3}-{7r\over4}} +\ \!V^{-{5r\over12}}\Big)V^{6\epsilon}\bigg)\ \!=\ \! \sum_{j=1}^{J-1}{(-z_0\log U)^j\over j!}\ \!+\ \! {U^{z_0-z_1}(z_0\log U)^J\over J!}\ \!+\ \!O_{\epsilon}\Big(V^{-{5r\over12}+6\epsilon}\Big)}\leqno(44)$$

since $\ \!2r\le z_0\le201r/100\ \ \!,\ \ \!z_0\le z_1\le 2z_0\ \ \!,\ \ \! T\le V^{r\over2}\le V^{1/200}\ .$

Also $\ J/2 -2\ \!=\ \! [z_0\log U+2]-2\ \!\le\ \!z_0\log U <\ \![(z_0\log U+2]-1 = J/2 -1\ .\ $

Using (3) with $\ \!u= -z_0\log U ,\ \!$ and both the bounds of lemma 1 (the latter with $\ \!j= J$) we then

obtain
$$\eqalign{&\quad\ \ \sum_{j=1}^{J-1}{(-z_0\log U)^j\over j!}\ \!+\ \! {U^{z_0-z_1}(z_0\log U)^J\over J!}\ \!\le\ {(-z_0\log U)^{J-1}\over2(J-1)!}\ \!+\ \!{U^{z_0-z_1}(z_0\log U)^J\over J!}\ \!\le\cr
&\qquad\qquad\ \ \le\ \! -{(z_0\log U)^{J-1}\over(J-1)!}\Big({1\over2}-{z_0\log U\over J}\Big)\ \!\le\ \!-{(z_0\log U)^{J-1}\over(J-1)!}\Big({1\over2}-{J-2\over2J}\Big)\ =\cr
&=\ \!-{(z_0\log U)^{J-1}\over J!}\ \le\ \!-{1\over J!}\Big({J-4\over2}\Big)^{J-1}\le\ \!-{1\over J^2}\Big({e(J-4)\over2J}\Big)^{J-1} =\ -{1\over J^2}\Big({e\over2}\Big)^{J-1}\Big({J-4\over J}\Big)^3\ \!\cdot\cr
&\qquad\quad\ \cdot\Big(1-{4\over J}\Big)^{J-4} \le\ -{1\over e^4J^2}\Big({J-4\over J}\Big)^3\Big({e\over2}\Big)^{J-1}\le\ -{1\over e^4J^2}\Big({J-4\over J}\Big)^3V^{4(1-\log2)z_0\over3}}\leqno(45)$$

since $J-1>2z_0\log U = (4z_0\log V)/3\ \!$ and since the sequence $\ \!\big((j-4)/j\big)^{j-4}$ is decreasing for

$j\ge5$ and has limit $e^{-4}\ .$ Therefore, by (44), (45)
$$V^{4(1-\log2)z_0\over3}\ \!\le\ e^4 J^2\Big({J\over J-4}\Big)^3\bigg(1\ \!+\ \!O_{\epsilon}\Big(V^{-{5r\over12}+6\epsilon}\Big)\bigg)\ .$$

Moreover $\ J^2\le 4\big(z_0\log U +2\big)^2 = 16\big(z_0\log V+3\big)^2/9\ $ while, according to lemma 5, we have

$J>2z_0\log U\ge 20\ \!,\ $ so that $\ J^3(J-4)^{-3}\le 125/64< 63/32\ .\ $ Hence
$${V^{4(1-\log2)z_0\over3}\over(z_0\log V+3)^2}\ \!\le\ {16e^4\over9}\Big({J\over J-4}\Big)^3\bigg(1\ \!+\ \!O_{\epsilon}\Big(V^{-{5r\over12}+6\epsilon}\Big)\bigg)\ \!\!<\  \!{7e^4\over 2}+\ \!O_{\epsilon}\Big(V^{-{5r\over12}+6\epsilon}\Big)\ .$$

But this is impossible when $V$ is large enough.

We then have either $b\ = 0\ $ or $\  b= 1/2\ .$

On the other hand QRH implies $\ \! b<1/2\ \!.\ \!$ Hence $\ \!b = 0\ \!$ and the Riemann Hypothesis is true.

\qed

\non

\centerline{\bf References}

\non

\item{[1]} E. Bombieri, {\it Le Grand Crible dans la Th\'eorie Analytique des Nombres} (seconde edition revue et augment\'ee), Ast\'erisque {\bf 18} (1974)

\item{[2]} E. C. Titchmarsh, {\it The theory of the Riemann Zeta--function} (second edition revised by D.R. Heath--Brown), Clarendon Press, Oxford, 1986

\end